\documentclass[12pt,leqno]{amsart}
\usepackage{latexsym, amssymb, amsmath, amscd, eucal}

\theoremstyle{plain}
    \newtheorem{rema}{Remark}[section]
    \newtheorem{propo}[rema]{Proposition}

   \newtheorem{theo}[rema]{Theorem}
 \newtheorem{conj}[rema]{Conjecture}
   \newtheorem{defi}[rema]{Definition}
    \newtheorem{lemma}[rema]{Lemma}
    \newtheorem{corol}[rema]{Corollary}
     \newtheorem{exam}[rema]{Example}
  \newtheorem{rmk}[rema]{Remark}

 \newtheorem{prob}[rema]{Open Problem}

	\newcommand{\nno}{\nonumber}

	\newcommand{\p}{\partial}

 \newcommand{\pf}{{\it Proof:}\hspace{2ex}}
 
 \newcommand{\epfv}{\hspace{1em}$\Box$\vspace{1em}}

\newcommand{\bZ}{{\mathbb Z}}

\newcommand{\bN}{{\mathbb N}}

\newcommand{\BQ}{\begin{eqnarray}}
\newcommand{\EQ}{\end{eqnarray}}
\newcommand{\BQn}{\begin{eqnarray*}}
\newcommand{\EQn}{\end{eqnarray*}}

\newcommand{\BL}{\begin{align}}
\newcommand{\EL}{\end{align}}
\newcommand{\BLn}{\begin{align*}}
\newcommand{\ELn}{\end{align*}}

\newcommand{\Hes}{ \text{Hes\,} }

\newcommand{\bC}{\mathbb C}
\newcommand{\bCzz}{\bC[[z]]}
\newcommand{\bCz}{\bC[z]}
\newcommand{\bT}{{\mathbb T}}

\newcommand{\cA}{{\mathcal A}}
\newcommand{\bG}{{\mathbb G}}
\newcommand{\cG}{{\mathcal G}}

\newcommand{\la}{\langle}
\newcommand{\ra}{\rangle}

\renewcommand{\theequation}{\thesection.\arabic{equation}}
\renewcommand{\therema}{\thesection.\arabic{rema}}
\newcommand{\EAn}{\end{align*}}

\title[Some Properties of and Open Problems on HNPs]
{Some Properties of and Open Problems on 
Hessian Nilpotent Polynomials}

    \author{Wenhua Zhao}      
    \date{\today}
    
\begin{document}

\textwidth=13.5cm
\baselineskip=16pt

\begin{abstract}

\textwidth=13.5cm
\baselineskip=16pt

In the recent work \cite{BE1}, \cite{Me}, \cite{Burgers} and \cite{HNP}, 
the well-known Jacobian conjecture (\cite{BCW}, \cite{E}) 
has been reduced to a problem on HN (Hessian nilpotent) 
polynomials (the polynomials whose Hessian matrix are nilpotent) 
and their (deformed) inversion pairs. 
In this paper, we prove several results on HN  
polynomials, their (deformed) inversion pairs as well as on
the associated symmetric polynomial or formal maps. 
We also propose some open problems for 
further study of these objects. 
\end{abstract}

   

\thanks{The author has been partially supported 
by NSA Grant R1-07-0053}

 \bibliographystyle{alpha}
    \maketitle


\renewcommand{\theequation}{\thesection.\arabic{equation}}
\renewcommand{\therema}{\thesection.\arabic{rema}}
\setcounter{equation}{0}
\setcounter{rema}{0}
\setcounter{section}{0}

\section{\bf Introduction}\label{S1}

In the recent work \cite{BE1}, \cite{Me}, \cite{Burgers} and \cite{HNP}, 
the well-known Jacobian conjecture (see \cite{BCW} and \cite{E}) 
has been reduced to a problem on HN (Hessian nilpotent) 
polynomials, i.e.\,the polynomials whose 
Hessian matrix are nilpotent, 
and their (deformed) inversion pairs. In this paper, 
we prove some properties of HN polynomials, 
the (deformed) inversion pairs of (HN) polynomial, 
the associated symmetric polynomial 
or formal maps, the graphs assigned 
to homogeneous harmonic polynomials, 
etc. Another purpose of this paper is 
to draw the reader's attention to some open problems 
which we believe will be interesting and important 
for further study of these objects.

In this section we first discuss
some backgrounds and motivations 
in Subsection \ref{S1.1} for 
the study of HN polynomials and their (deformed) 
inversion pairs. We also fix some 
terminology and notation in this subsection 
that will be used throughout this paper. 
Then in Subsection \ref{S1.2} we give an  
arrangement description of this paper.

\subsection{Background and Motivation}\label{S1.1}

Let $z=(z_1, z_2, \dots, z_n)$ be 
$n$ free commutative variables. We denote 
by $\bCz$ (resp.\,$\bCzz$) the 
algebra of polynomials (resp.\,formal power series)
of $z$ over $\bC$. A polynomial or 
formal power series $P(z)$ is said to be 
{\it HN} ({\it Hessian nilpotent}) 
if its Hessian matrix 
$\Hes P\!:=(\frac{\p^2 P}{\p z_i \p z_j})$ 
are nilpotent. The study of HN polynomials 
is mainly motivated by the recent progress achieved in
\cite{BE1}, \cite{Me}, \cite{Burgers} and \cite{HNP} on the well-known 
JC (Jacobian conjecture), which we will 
briefly explain below.

Recall that the JC first proposed 
by Keller \cite{Ke} in 1939 
claims: {\it for any polynomial map $F$ 
of $\bC^n$ with the Jacobian $j(F)=1$, 
its formal inverse map $G$ must also be a polynomial map}. 
Despite intense study 
for more than half a century, 
the conjecture is still open
even for the case $n=2$. 
For more history and known results 
before $2000$ on the Jacobian conjecture, 
see \cite{BCW}, \cite{E} 
and references there.
In $2003$, M. de Bondt, A. van den Essen 
(\cite{BE1}) and G. Meng (\cite{Me}) independently made 
the following breakthrough on the JC.

Let $D_i\!:=\frac{\p}{\p z_i}$ $(1\leq i\leq n)$ 
and $D=(D_1,  \dots, D_n)$. For any $P(z)\in \bCzz$, 
denote by $\nabla P(z)$ the {\it gradient} of $P(z)$, 
i.e. $\nabla P(z):=(D_1 P(z), \dots, D_n P(z))$.
We say a formal map $F(z)=z-H(z)$ is 
{\it symmetric} if $H(z)=\nabla P(z)$ 
for some $P(z)\in \bCzz$. Then, 
the {\it symmetric reduction} of the JC 
achieved in \cite{BE1} and \cite{Me} 
is that, {\it to prove or disprove 
the JC, it will be enough to consider 
only symmetric polynomial maps}. Combining with 
the classical {\it homogeneous reduction}
achieved in \cite{BCW} and \cite{Y},  
one may further assume that 
{\it the symmetric polynomial maps 
have the form $F(z)=z-\nabla P(z)$ with 
$P(z)$ homogeneous $($of degree $4$$)$.}
Note that, in this case the Jacobian condition 
$j(F)=1$ is equivalent to the condition that 
$P(z)$ is HN. For some other recent results 
on symmetric polynomial or formal maps, 
see \cite{BE1}--\cite{BE5}, 
\cite{EW},  \cite{Me}, \cite{Wr1}, 
\cite{TreeExp}, \cite{Burgers}, 
\cite{HNP} and \cite{EZ}.

Based on the homogeneous reduction and the symmetric reduction
of the JC discussed above, the author further showed in \cite{HNP} that
the JC is actually equivalent to the following so-called
{\it vanishing conjecture} of HN polynomials.

\begin{conj} \label{VC} $(${\bf Vanishing Conjecture}$)$ 
Let $\Delta\!:=\sum_{i=1}^n D_i^2$ be the Laplace operator of $\bCz$.
Then, for any HN polynomial $P(z)$ $($of homogeneous of degree $d=4$$)$, 
$\Delta^m P^{m+1}(z)=0$ when $m>>0$.
\end{conj}

Furthermore, the following criterion of Hessian nilpotency for 
formal power series was also proved in \cite{HNP}.

\begin{propo}\label{Crit-1}
For any $P(z)\in \bC[[z]]$ with $o(P(z))\geq 2$, 
the following statements are equivalent.
\begin{enumerate}
\item[(1)] $P(z)$ is HN.
\item[(2)] $\Delta^m P^m=0$ for any $m\geq 1$.
\item[(3)] $\Delta^m P^m=0$ for any $1\leq m\leq n$.
\end{enumerate}
\end{propo}

One crucial idea of the proofs in \cite{HNP} 
for the results above is to study a special 
formal deformation of symmetric formal maps. 
More precisely, let $t$ be a central 
formal parameter. For any $P(z)\in \bCzz$, 
we call $F(z)=z-\nabla P(z)$ 
the {\it associated symmetric maps} of $P(z)$. 
Let $F_t(z)\:=z-t\nabla P(z)$.
When the order $o(P(z))$ of $P(z)$ with respect to $z$ 
is greater than or equal to $2$, $F_t(z)$ is a formal map 
of $\bC[[t]][[z]]$ with $F_{t=1}(z)=F(z)$. Therefore, 
we may view $F_t(z)$ as a formal deformation 
of the formal map $F(z)$. In this case, 
one can also show (see \cite{Me} or 
Lemma $3.14$ in \cite{Burgers}) that 
the formal inverse map $G_t(z)\!:=F_t^{-1}(z)$ 
of $F_t(z)$ does exist and is also symmetric, 
i.e. there exists a unique $Q_t(z)\in \bC[[t]][[z]]$ 
with $o(Q_t(z))\geq 2$ such that 
$G_t(z)=z+t\nabla Q_t(z)$. We call $Q_t(z)$ 
the {\it deformed inversion pair} of $P(z)$. Note that, 
whenever $Q_{t=1}(z)$ makes sense, the formal inverse map $G(z)$ 
of $F(z)$ is given by $G(z)=G_{t=1}(z)=z+ \nabla Q_{t=1}(z)$, 
so in this case we call $Q(z) \!:=Q_{t=1}(z)$ the 
{\it inversion pair}\label{Inv-pair} of $P(z)$.

Note that, under the condition $o(P(z))\ge 2$, 
the deformed inversion pair $Q_t(z)$ of $P(z)$ 
might not be in $\bC[t][[z]]$, so $Q_{t=1}(z)$ 
may not make sense. But, if we assume further that 
$J(F_t)(0)=1$, or equivalently, $(\Hes P)(0)$ 
is nilpotent, then $F_t(z)$ is an automorphism 
of $\bC[t][[z]]$, hence so is its inverse map $G_t(z)$. 
Therefore, in this case $Q_t(z)$ lies in $\bC[t][[z]]$ 
and $Q_{t=1}(z)$ makes sense. Throughout this paper, 
whenever the inversion pair $Q(z)$ of a polynomial 
or formal power series $P(z)\in \bC[[z]]$ (not necessarily HN)
is under concern, our assumption on $P(z)$ will always be 
$o(P(z))\ge 2$ and $(\Hes P)(0)$ is nilpotent. Note that, 
for any HN $P(z)\in \bC[[z]]$ with $o(P(z))\ge 2$, 
the condition that $(\Hes P)(0)$ is nilpotent holds 
automatically.

For later purpose, let us recall the following 
formula derived in \cite{HNP} 
for the deformed inversion pairs 
of HN formal power series.

\begin{theo}\label{E-Q}
Suppose $P(z)\in \bC[[z]]$ with $o(P(z))\geq 2$ is HN. Then, 
 we have
\begin{align}\label{E-Q-1}
Q_t (z)&=\sum_{m=0}^\infty \frac {t^m}{2^m m!(m+1)!} \Delta^m P^{m+1}(z), 
\end{align}
\end{theo}

From the equivalence of the JC and the VC 
discussed above, we see that the study on 
the HN polynomials and their (deformed) inversion pairs 
becomes important and necessary, 
at least when the JC is concerned. 
Note that, due to the identity 
$\mbox{Tr\,}\Hes P = \Delta P$, 
HN polynomials are just a special family of
harmonic polynomials which are among 
the most classical objects 
in mathematics. Even though harmonic polynomials 
had been very well studied since the 
late of the eighteenth century, it seems that 
not much has been known on HN polynomials. 
We believe that these mysterious (HN) polynomials 
deserve much more attentions from mathematicians.

\subsection{Arrangement}\label{S1.2}
Considering the length of this paper, we here give 
a more detailed arrangement description 
of the paper. 

In Section \ref{dis-cncd}, 
we consider the following 
two questions. Let $P, S, T \in \bCzz$ 
with $P=S+T$ and $Q, U, V$ 
their inversion pairs, respectively. 

${\bf Q_1}$: 
{\it Under what conditions,  $P$ is HN iff 
both $S$ and $T$ are HN?} 

${\bf Q_2}$: 
{\it Under what conditions, we have $Q=U+V$?}

We give some sufficient conditions in Theorems 
\ref{p-inv-MainThm-1} and 
\ref{p-inv-MainThm-2} for the 
two questions above.
In Section \ref{cvg}, we employ a recursion formula
of inversion pairs derived in \cite{Burgers}
and Eq.\,$(\ref{E-Q-1})$ above to derive 
some estimates for the radius of 
convergence of inversion pairs 
of homogeneous (HN) polynomials 
(see Propositions \ref{cvg-MainThm-1} and 
\ref{cvg-MainThm-2}).

For any $P(z)\in \bCzz$, we say it is 
{\it self-inverting} if its inversion pair $Q(z)$ 
is $P(z)$ itself. In Section \ref{slf}, 
by using a general result on quasi-translations 
proved in \cite{B2}, we derive 
some properties of HN {\it self-inverting} 
formal power series $P(z)$. Another purpose of 
this section is to draw the reader's attention to 
Open Problem \ref{slf-prob} on classification of 
HN self-inverting polynomials 
or formal power series. 

In Section \ref{pVC}, we show in 
Proposition \ref{pVC-MainThm}, 
when the base field has characteristic 
$p>0$, the VC, unlike the JC, 
actually holds for any polynomials $P(z)$ 
even without the HN condition on $P(z)$.
It also holds in this case for any HN formal 
power series.
One interesting question 
(see Open Problem \ref{pVC-prob}) 
is to see if the VC like the JC fails over $\bC$ 
when $P(z)$ is allowed to be any 
HN formal power series.

In Section \ref{crn}, we prove a criterion 
of Hessian nilpotency for  
homogeneous polynomials over $\bC$ 
(see Theorem \ref{crn-MainThm}). 
Considering the criterion in Proposition \ref{Crit-1},
this criterion is somewhat surprising but its proof 
turns out to be very simple.

Section \ref{symmc} is mainly motivated 
by the following question raised 
by M. Kumar (\cite{Ku}) and 
D. Wright (\cite{Wr3}). Namely, 
for a symmetric formal map 
$F(z)=z-\nabla P(z)$, how to write 
$f(z)\!:=\frac 12 \sigma_2 -P(z)$ 
(where $\sigma_2\!:=\sum_{i=1}^n z_i^2$)
and $P(z)$ itself as formal 
power series in $F(z)$? 
In this section, we derive some explicit formulas 
to answer the questions above and also for the 
same question for $\sigma_2$ 
(see  Proposition \ref{symmc-MainThm-1}). 
From these formulas, 
we also show in Theorem \ref{symmc-MainThm-2}
that, the VC holds for a HN polynomial $P(z)$ 
iff one (hence, all) of $\sigma_2$, $P(z)$ 
and $f(z)$ can be written as a polynomial in $F$, 
where $F(z)=z-\nabla P(z)$ is the associated 
polynomial maps of $P(z)$.  

Finally, in Section \ref{grf}, we discuss 
a graph $\cG(P)$ assigned to 
each homogeneous harmonic polynomials $P(z)$. 
The graph $\cG(P)$ was first proposed by the author 
and later was further studied by Roel Willems  
in his master thesis \cite{Roel} 
under direction of Professor Arno van den Essen.
In Subsection \ref{S8.1} we give the definition 
of the graph $\cG(P)$ for any homogeneous  
harmonic polynomial $P(z)$ and discuss 
the {\it connectedness reduction} 
(see Corollary \ref{conn-reduction})
which says, to study the VC for homogeneous HN
polynomials $P(z)$, it will be enough to consider the case
when the graph $\cG(P)$ is connected.
In Subsection \ref{S8.2} we consider a connection of $\cG(P)$ 
with the tree expansion formula derived 
in \cite{Me} and \cite{TreeExp} for the inversion 
pair $Q(z)$ of $P(z)$ (see also Proposition 
\ref{grf-MainThm-2}). As an application of the connection, 
we use it to give another proof for 
the connectedness reduction discussed in 
Corollary \ref{conn-reduction}.

One final remark on the paper is as follows. 
Even though we could have focused only 
on (HN) polynomials, at least when only the JC is concerned, 
we will formulate and prove our results 
in the more general setting of (HN) 
formal power series whenever 
it is possible. \\

{\bf Acknowledgement:} The author is very grateful to 
Professors Arno van den Essen, Mohan Kumar and David Wright 
for inspiring communications and constant encouragement. 
Section \ref{symmc} was mainly motivated 
by some questions raised by Professors 
Mohan Kumar and David Wright. The author also would like 
to thank Roel Willems for sending the author his master 
thesis in which he has obtained some very interesting 
results on the graphs $\cG(P)$ of homogeneous 
harmonic polynomials. At last but not the least, 
the author thanks the referee and the editor for 
many valuable suggestions.


\renewcommand{\theequation}{\thesection.\arabic{equation}}
\renewcommand{\therema}{\thesection.\arabic{rema}}
\setcounter{equation}{0}
\setcounter{rema}{0}

\section{\bf Disjoint  Formal Power Series 
and Their Deformed Inversion Pairs}
\label{dis-cncd}

Let $P, S, T\in \bCzz$ with $P=S+T$, and 
$Q$, $U$ and $V$ their inversion pairs, 
respectively. In this section, 
we consider the following two questions: 

\begin{enumerate}
\item[${\bf{Q_1}}$:] \label{Q1}
{\it Under what conditions, 
 $P$ is HN if and only if both $S$ and $T$ are HN?}

\item[${\bf Q_2}$:] \label{Q2}
{\it  Under what conditions, we have $Q=U+V$?}
\end{enumerate}

We give some answers to the questions 
${\bf Q_1}$ and ${\bf Q_2}$ in 
Theorems \ref{p-inv-MainThm-1} 
and \ref{p-inv-MainThm-2}, respectively.
The results proved here will also be needed 
in Section \ref{grf} when we consider  
a graph associated to 
homogeneous harmonic polynomials. 

To question ${\bf Q_1}$ above, 
we have the following result.

\begin{theo}\label{p-inv-MainThm-1}
Let $S, T\in \bCzz$ such that  
$\la\nabla (D_i S), \nabla (D_j T)\ra=0$ for any $1\leq i, j\leq n$, 
where $\la\cdot, \cdot\ra$ denotes the standard 
$\bC$-bilinear form of $\bC^n$. 
Let $P=S+T$. Then, we have  
\begin{enumerate}
\item[$(a)$] $\Hes(S)\, \Hes(T)=\Hes(T)\, \Hes(S)=0.$

\item[$(b)$] $P$ is HN iff both $S$ and $T$ are HN.
\end{enumerate}
\end{theo}

Note that statement $(b)$ in the theorem above 
was first proved by R. Willems  (\cite{Roel}) 
in a special setting as in Lemma \ref{p-inv-L4} below  
for homogeneous harmonic polynomials.

\pf $(a)$ For any $1\leq i, j\leq n$, consider the $(i, j)^{th}$ 
entry of the product $\Hes (S) \Hes(T)$:
\begin{align}
\sum_{k=1}^n \frac {\p^2 S}{\p z_i \p z_k} \frac {\p^2 T}{\p z_k \p z_j}
=\la \nabla (D_i S), \nabla (D_j T)\ra =0.
\end{align}

Hence $\Hes (S)\, \Hes(T)=0$. Similarly, we have 
$\Hes(T)\, \Hes(S)=0$.

$(b)$ follows directly from $(a)$ and the lemma below.
\epfv

\begin{lemma}\label{p-inv-L1}
Let $A$, $B$ and $C$ be $n\times n$ matrices 
with entries in any commutative ring. Suppose that 
$A=B+C$ and $BC=CB=0$. Then, $A$ 
is nilpotent iff both $B$ and $C$ 
are nilpotent.
\end{lemma}

\pf The $(\Leftarrow)$ part is trivial 
because $B$ and $C$ in particular commute 
with each other. 

To show $(\Rightarrow)$, note that $BC=CB=0$. So for any $m\geq 1$, 
we have
\begin{align*}
A^m B=(B+C)^m B=(B^m+C^m)B=B^{m+1}.
\end{align*}
Similarly, we have $C^{m+1}=A^mC$. Therefore, 
if $A^N=0$ for some $N\geq 1$, 
we have $B^{N+1}=C^{N+1}=0$.
\epfv

Note that, for the $(\Leftarrow)$ part of $(b)$ 
in Theorem \ref{p-inv-MainThm-1}, we need only 
a weaker condition. Namely, 
for any $1\leq i, j\leq n$, 
\begin{align*}
\la\nabla (D_i S), \nabla (D_j T)\ra =
\la\nabla (D_j S), \nabla (D_i T)\ra,
\end{align*}
which will ensure that $\Hes (S)$ and $\Hes (T)$
commute.

To consider the second question ${\bf Q_2}$, 
let us first fix the following notation.

For any $P\in \bCzz$, let $\cA(P)$ denote 
the subalgebra  of $\bCzz$ generated 
by all partial derivatives of $P$ 
(of any order). We also define 
a sequence $\{ Q_{[m]}(z) \,|, m\geq 1\}$ 
by writing the deformed inversion pair 
$Q_t(z)$ of $P(z)$ as
\begin{align}\label{Qm}
Q_t(z)=\sum_{m\geq 1} t^{m-1} Q_{[m]}(z).
\end{align}

\begin{lemma}\label{p-inv-L2}
For any $P\in \bCzz$, we have 

$(a)$ $\cA(P)$ is closed 
under the action of any differential operator 
of $\bCz$ with constant coefficients.

$(b)$ For any $m\ge 1$, we have $Q_{[m]}(z)\in \cA (P)$. 
\end{lemma}

\pf $(a)$ Note that, by the definition of $\cA(P)$,
a formal power series $g(z)\in \bCzz$ 
lies in $\cA(P)$ iff it can be written 
(not necessarily uniquely) as a polynomial 
in partial derivatives of $P(z)$. 
Then, by the Leibniz Rule, it is easy 
to see that, for any $g(z)\in \cA(P)$, 
$D_i g(z)\in \cA(P)$ $(1\leq i \leq n)$.
Repeating this argument, we see that 
any partial derivative of $g(z)$ is 
in $\cA(P)$. Hence $(a)$ follows.

$(b)$ Recall that, by Proposition $3.7$ 
in \cite{Burgers},
we have the following 
recurrent formula for 
$Q_{[m]}(z)$ $(m\geq 1)$ 
in general:
\allowdisplaybreaks{
\begin{align}
Q_{[1]}(z)&=P(z),\label{Recurr-1} \\
Q_{[m]}(z)&=\frac 1{2(m-1)} 
\sum_{\substack {k, l\geq 1 \\ k+l=m}}
\la \nabla Q_{[k]}(z), \nabla Q_{[l]}(z)\ra. 
\label{Recurr-2}
\end{align} }
for any $m\geq 2$.

By using $(a)$, the recurrent formulas above and 
induction on $m\geq 1$, it is easy to check  
that $(b)$ holds too.
\epfv

\begin{defi} For any $S, T\in \bCzz$, we say 
$S$ and $T$ are {\it disjoint} to each other 
if, for any $g_1\in \cA(S)$ and $g_2\in \cA(T)$, 
we have $\la \nabla g_1, \nabla g_2 \ra =0$. 
\end{defi}

This terminology will be justified in 
Section \ref{grf} when we consider 
a graph $\cG(P)$ associated to 
homogeneous harmonic polynomials $P$.
 
\begin{lemma}\label{p-inv-L3}
Let $S, T\in \bCzz$. Then $S$ and $T$ are disjoint to each other iff, 
for any $\alpha, \beta \in \bN^n$, we have
\begin{align}\label{p-inv-L3-e1}
\la \nabla (D^\alpha S), 
\nabla (D^\beta T) \ra=0. 
\end{align}
\end{lemma}

\pf The $(\Rightarrow)$ part of the lemma is trivial. 
Conversely, for any $g_1\in \cA(S)$ and $g_2\in \cA(T)$
$(i=1, 2)$, we need show
\begin{align*}
\la \nabla g_1, \nabla g_2 \ra=0. 
\end{align*}
But this can be easily checked by, first, 
reducing to the case that $g_1$ and $g_2$ 
are monomials of partial derivatives 
of $S$ and $T$, respectively, 
and then applying the Leibniz rule 
and Eq.\,$(\ref{p-inv-L3-e1})$ above. 
\epfv

A family of examples of disjoint 
polynomials or formal power series 
are given as in the following lemma, 
which will also be needed later in 
Section \ref{grf}. 

\begin{lemma}\label{p-inv-L4}
Let $I_1$ and $I_2$ be two finite subsets 
of $\bC^n$ such that, for any 
$\alpha_i \in I_i$ $(i=1, 2)$, 
we have $\la \alpha_1, \alpha_2 \ra=0$.
Denote by $\cA_i$ $(i=1, 2)$ the completion 
of the subalgebra of $\bCzz$ generated by 
$h_\alpha (z)\!:=\la \alpha, z \ra$ 
$(\alpha\in I_i)$, i.e. 
$\cA_i$ is the set of all formal power series in 
$h_\alpha (z)$ $(\alpha\in I_i)$ over $\bC$. 
Then, for any $P_i \in \cA_i$ $(i=1, 2)$, 
$P_1$ and $P_2$ are disjoint.
\end{lemma}

\pf First, by a similar argument as 
the proof for Lemma \ref{p-inv-L2}, 
$(a)$, it is easy to check that $\cA_i$ 
$(i=1, 2)$ are closed under 
action of any differential 
operator with constant 
coefficients. Secondly, 
since  $\cA_i$ 
$(i=1, 2)$ are subalgebras 
of $\bCzz$, we have $\cA(P_i)
\subset \cA_i$ $(i=1, 2)$. 

Therefore,  to show 
$P_1$ and $P_2$ are disjoint to each other, 
it will be enough to show that, 
for any $g_i \in \cA_i$ $(i=1, 2)$, 
we have $\la \nabla g_1, \nabla g_2\ra=0$.
But this can be easily checked 
by first reducing to the case when 
$g_i$ $(i=1, 2)$ are monomials 
of $h_\alpha(z)$ $(\alpha\in I_i)$, 
and then applying the Leibniz rule 
and the following identity: for any 
$\alpha, \beta\in \bC^n$, 
\begin{align*}
\la \nabla h_\alpha(z), \nabla h_\beta(z)\ra = \la \alpha, \beta\ra.
\end{align*}
\epfv

Now, for the second question ${\bf Q_2}$ 
on page \pageref{Q2}, we have the following result. 

\begin{theo}\label{p-inv-MainThm-2}
Let $P, S, T\in \bCzz$ with order greater than or equal to $2$, 
and $Q_t, U_t, V_t$ their deformed 
inversion pairs, respectively. 
Assume that $P=S+T$ and 
$S$, $T$ are disjoint 
to each other. 
Then 

$(a)$ $U_t$ and $V_t$ are also disjoint to each other, i.e. 
for any $\alpha, \beta\in \bN^n$, we have 
\begin{align*}
\left \la \nabla D^\alpha U_t(z), \nabla D^\beta V_t (z)\right \ra = 0.
\end{align*}

$(b)$ We further have 
\begin{align}
Q_t=U_t+V_t.
\end{align}
\end{theo}

\pf $(a)$ follows directly from 
Lemma \ref{p-inv-L2}, $(b)$ and Lemma \ref{p-inv-L3}.

$(b)$ Let $Q_{[m]}$, $U_{[m]}$ and 
$V_{[m]}$ $(m\geq 1)$ be defined 
as in Eq.\,$(\ref{Qm})$. Hence 
it will be enough to show 
\begin{align}\label{p-inv-MainThm-e1}
Q_{[m]}=U_{[m]}+V_{[m]}
\end{align}
for any $m\geq 1$.

We use induction on $m\geq 1$. 
When $m=1$, Eq.\,$(\ref{p-inv-MainThm-e1})$ 
follows from the condition $P=S+T$ and 
Eq.\,$(\ref{Recurr-1})$ . 
For any $m\geq 2$, by Eq.\,$(\ref{Recurr-2})$  
and the induction assumption, we have
\allowdisplaybreaks{
\begin{align*}
 Q_{[m]}&=\frac 1{2(m-1)} 
\sum_{\substack {k, l\geq 1 \\ k+l=m}}
\la \nabla Q_{[k]}, \nabla Q_{[l]} \ra \\
&=\frac 1{2(m-1)} 
\sum_{\substack {k, l\geq 1 \\ k+l=m}}
\la \nabla U_{[k]} + \nabla V_{[k]}, 
\nabla U_{[l]} + \nabla V_{[l]} \ra\\
\intertext{Noting that, by Lemma \ref{p-inv-L2},  
$U_{[j]} \in \cA(S)$ and $V_{[j]} \in \cA(T)$ $(1\leq j\leq m)$:} 
&=\frac 1{2(m-1)} 
\sum_{\substack {k, l\geq 1 \\ k+l=m}}
\la \nabla U_{[k]}, 
\nabla U_{[l]} \ra 
+\frac 1{2(m-1)} 
\sum_{\substack {k, l\geq 1 \\ k+l=m}}
\la  \nabla V_{[k]}, 
\nabla V_{[l]} \ra\\
\intertext{Applying the recursion formula 
Eq.\,$(\ref{Recurr-2})$  to both 
$U_{[m]}$ and $V_{[m]}$:} 
&=U_{[m]}+V_{[m]}.
\end{align*} }
\epfv

As later will be pointed out in Remark \ref{PST-2},
one can also prove this theorem by using a tree 
expansion formula of inversion pairs, which was 
derived in \cite{Me} and \cite{TreeExp}, in the setting 
as in Lemma \ref{p-inv-L4}.

From Theorems \ref{p-inv-MainThm-1},  
\ref{p-inv-MainThm-2} and 
Eqs.\,$(\ref{E-Q-1})$, $(\ref{Qm})$, 
it is easy to see that 
we have the following corollary.

\begin{corol}\label{p-inv-MainCorol}
Let $P_i\in \bCzz$ $(1\leq i\leq k)$ 
which are disjoint to each other.
Set $P=\sum_{i=1}^k P_i$. Then, we have

$(a)$ $P$ is HN iff each $P_i$ is HN.

$(b)$ Suppose that $P$ is HN. Then, 
for any $m\geq 0$, 
we have
\begin{align}
\Delta^m P^{m+1}=\sum_{i=1}^k \Delta^m P_i^{m+1}.
\end{align}
Consequently, if the VC holds 
for each $P_i$, 
then it also holds for $P$.  
\end{corol}

\renewcommand{\theequation}{\thesection.\arabic{equation}}
\renewcommand{\therema}{\thesection.\arabic{rema}}
\setcounter{equation}{0}
\setcounter{rema}{0}

\section{\bf Local Convergence of Deformed Inversion Pairs 
of Homogeneous (HN) Polynomials}\label{cvg}

Let $P(z)$ be a formal power series 
which is convergent near $0\in \bC^n$.
Then the associated symmetric map $F(z)=z-\nabla P$ is 
a well-defined analytic map from an open neighborhood 
of $0\in \bC^n$ to $\bC^n$. 
If  we further assume that $JF(0)=I_{n\times n}$,  
the formal inverse $G(z)=z+\nabla Q(z)$ of $F(z)$ is 
also locally well-defined analytic map. 
So the inversion pair $Q(z)$ of $P(z)$ is also locally 
convergent near $0\in \bC^n$. 
In this section, we use the formulas Eqs.\,$(\ref{Recurr-2})$, 
$(\ref{E-Q-1})$ and the Cauchy estimates 
to derive some estimates for the radius of convergence 
of inversion pairs $Q(z)$ of homogeneous (HN) 
polynomials $P(z)$ (see Propositions \ref{cvg-MainThm-1} 
and \ref{cvg-MainThm-2}). 

First let us fix the following notation.

For any $a\in \bC^n$ and $r>0$, we denote by 
$B(a, r)$ (resp.\,$S(a, r)$) the open ball (resp.\,the sphere)
centered at $a\in \bC$ with radius $r>0$. The unit sphere 
$S(0, 1)$ will also be denoted by $S^{2n-1}$.
Furthermore, we let  $\Omega(a, r)$ 
be the polydisk centered at 
$a\in \bC^n$ with radius $r>0$, i.e.
$\Omega(a, r)\!:=\{ z\in \bC^n \,| \, 
|z_i-a_i|<r, \, 1\leq i\leq n \}$. 
For any subset $A \subset \bC^n$, we will use 
$\bar{A}$ to denote the closure of $A$ 
in $\bC^n$.

For any polynomial $P(z)\in \bCz$ and 
a compact subset $D\subset \bC^n$, 
we set $|P|_D$ to be the maximum value of $|P(z)|$ over $D$.
In particular, when $D$ is the unit sphere $S^{2n-1}$, 
we also write $|P|=|P|_D$, i.e.
\begin{align}\label{Def-|P|}
|P|\!:=\max \{|P(z)|\,|\, z\in S^{2n-1} \}.
\end{align}

Note that, for any $r\geq 0$ and $a\in B(0, r)$, we have
$\Omega(a, r) \subset B(a, r)\subset B(0, 2r)$.
Combining with the well-known Maximum Principle 
of holomorphic functions, we get
\begin{align}\label{cvg-E1}
|P|_{\overline{\Omega(a, r)}}
\leq |P|_{\overline{B(a, r)}}
\leq |P|_{\overline{B(0, 2r)}}
=|P|_{S(0, 2r)}.
\end{align}
 
For the inversion pairs $Q$ of homogeneous polynomials $P$
without HN condition, we have the following estimate 
for the radius of convergence at $0\in \bC^n$.

\begin{propo}\label{cvg-MainThm-1}
Let $P(z)$ be a non-zero 
homogeneous polynomial 
$($not necessarily HN$)$ 
of degree $d\geq 3$ 
and $r_0\:=(n 2^{d-1} |P|)^{\frac 1{2-d}}$. 
Then the inversion pair 
$Q(z)$ converges over 
the open ball $B(0, r_0)$.
\end{propo}

To prove the proposition, we need the following lemma.

\begin{lemma}\label{cvg-LL1} 
Let $P(z)$ be any polynomial and $r>0$. Then, 
for any $a\in B(0, r)$ and $m\geq 1$, we have
\begin{align}\label{cvg-LL1-e}
\left |Q_{[m]}(a)\right | \leq \frac {n^{m-1} 
|P|^m_{S(0, 2r)}}{2^{m-1} r^{2m-2}}. 
\end{align}
\end{lemma}
\pf We use induction on $m\geq 1$. 
First, when $m=1$, by Eq.\,$(\ref{Recurr-1})$ we have $Q_{[1]}=P$. 
Then Eq.\,$(\ref{cvg-LL1-e})$ follows from 
the fact $B(a, r)\subset B(0, 2r)$ and the maximum 
principle of holomorphic functions. 

Assume Eq.\,$(\ref{cvg-LL1-e})$ holds for any $1\leq k\leq m-1$.
Then, by the Cauchy estimates of holomorphic functions 
(e.g. see Theorem $1.6$ in \cite{R}), 
we have
\begin{align}\label{cvg-LL1-pe-1}
\left |(D_i Q_{[k]}) (a)\right | 
\leq \frac 1r \left | Q_{[k]} \right |_{\overline{\Omega(0, r)}}
\leq \frac {n^{k-1} |P|^k_{B(0, 2r)}}{2^{k-1} r^{2k-1}}. 
\end{align}

By Eqs.\,$(\ref{Recurr-2})$ and  $(\ref{cvg-LL1-pe-1})$, we have
\allowdisplaybreaks{
\begin{align*}
|Q_{[m]}(a)|& \leq \frac 1{2(m-1)} 
\sum_{\substack {k, l\geq 1 \\ k+l=m}}
\left| \la \nabla Q_{[k]}, \nabla Q_{[l]} \ra \right | \\
&\leq \frac 1{2(m-1)} 
\sum_{\substack {k, l\geq 1 \\ k+l=m}}
n \, \frac {n^{k-1} |P|^k_{S(0, 2r)}}{2^{k-1} r^{2k-1}} 
\frac {n^{\ell-1} |P|^\ell_{S(0, 2r)}}{2^{\ell-1} r^{2\ell-1}} \\
&=\frac {n^{m-1} |P|^m_{S(0, 2r)}}{2^{m-1} r^{2m-2}}. 
\end{align*}
}\epfv

\underline{\it Proof of Proposition 
\ref{cvg-MainThm-1}}:
By Eq.\,$(\ref{Qm})$ , we know that, 
\begin{align}\label{cvg-MainThm-1-pe}
Q(z)=\sum_{m\geq 1} Q_{[m]}(z).
\end{align}

To show the proposition, it will be enough to show 
the infinite series above converges absolutely
over $B(0, r)$ for any $r<r_0$.
 
First, for any $m\geq 1$, let $A_m$ be the RHS of the 
inequality Eq.\,$(\ref{cvg-LL1-e})$. 
Note that, since $P$ is homogeneous of degree $d\geq 3$, 
we further have 
\begin{align}
|P|^m_{B(0, 2r)}=
\left ((2r)^d |P|_{S^{2n-1}} \right)^m =(2r)^{dm}|P|^m.
\end{align}
Therefore, for any $m\geq 1$, we have
\begin{align}
A_m=2^{(d-1)m +1} n^{m-1}r^{(d-2)m+2} |P|^m,
\end{align}
and by Lemma \ref{cvg-LL1}, 
\begin{align}
|Q_{[m]}(a)| \leq A_m
\end{align}
for any $a\in B(0, r)$.

Since $0<r<r_0=(n 2^{d-1}  |P|)^{2-d}$, 
it is easy to see that
\begin{align*}
\lim_{m\to +\infty} \frac{A_{m+1}}{A_m} = n2^{d-1}  r^{d-2}|P| < 1. 
\end{align*}

Therefore, by the comparison test, the infinite series
in Eq.\,$(\ref{cvg-MainThm-1-pe})$ 
converges absolutely and uniformly 
over the open ball $B(0, r)$. 
\epfv

Note that the estimate given in Proposition \ref{cvg-MainThm-1} 
depends on the number $n$ of variables. Next we show that, 
with the HN condition on $P$, an estimate independent  
of $n$ can be obtained as follows.

\begin{propo}\label{cvg-MainThm-2}
Let $P(z)$ be a homogeneous HN polynomial 
of degree $d\geq 4$ and set $r_0\!:=(2^{d+1}|P|)^{\frac 1{2-d}}$.
Then, the inversion pair $Q(z)$ of $P(z)$ 
converges over the open ball $B(0, r_0)$.  
\end{propo}

Note that, when $d=2$ or $3$, by Wang's Theorem (\cite{Wa}), 
the JC holds in general. Hence it also holds 
for the associated symmetric map $F(z)=z-\nabla P$ when $P(z)$ is HN. 
Therefore $Q(z)$ in this case is also a polynomial of $z$
and converges over the whole space $\bC^n$.

To prove the proposition above, we first need the following two lemmas.

\begin{lemma}\label{cvg-L1}
Let $P(z)$ be a homogeneous  
polynomial of degree $d\geq 1$ and $r>0$. 
For any $a \in B(0, r)$, $m\geq 0$ and $\alpha \in \bN^n$, 
we have 
\begin{align}\label{cvg-L1-e1}
|(D^\alpha P^{m+1})(a)| \leq \frac{\alpha!}{r^{ |\alpha|} } 
(2r)^{d(m+1)} |P|^{m+1}.
\end{align}
\end{lemma}
\pf
First, by the Cauchy estimates and Eq.\,$(\ref{cvg-E1})$, 
we have
\begin{align}\label{cvg-L1-pe1}
|(D^\alpha P^{m+1})(a)|& \leq \frac{\alpha!}{r^{ |\alpha|} }
|P^{m+1}|_{\overline{\Omega (a, r)}} \leq \frac{\alpha!}{r^{ |\alpha|} } 
|P^{m+1}|_{\overline{B(0, 2r)}}. 
\end{align}

On the other hand, by the maximum principle 
and the condition that $P$ is homogeneous of degree $d\geq 3$, we have 
\allowdisplaybreaks{
\begin{align}\label{cvg-L1-pe2}
|P^{m+1}|_{\overline{B(0, 2r)}} & = |P|_{\overline{B(0, 2r)} }^{m+1} 
=|P|_{\overline{S(0, 2r)} }^{m+1} =((2r)^d |P|)^{m+1}  \\
&=(2r)^{d(m+1)}|P|^{m+1}.\nno
\end{align} }

Then, combining Eqs.\,$(\ref{cvg-L1-pe1})$ and 
$(\ref{cvg-L1-pe2})$, we get 
Eq.\,$(\ref{cvg-L1-e1})$.
\epfv

\begin{lemma}\label{cvg-L2}
For any $m\geq 1$, we have
\begin{align}\label{cvg-L2-e1}
\sum_{\substack{\alpha \in \bN^n \\ |\alpha|=m }} \alpha ! 
\leq m!\binom{m+n-1}{m}=\frac{(m+n-1)!}{(n-1)!}.
\end{align}
\end{lemma}
\pf First, for any $\alpha\in \bN^n$ with $|\alpha|=m$, we have
$\alpha!\leq m!$ since the binomial $\binom{m}{\alpha}=\frac{m!}{\alpha!}$ 
is always a positive integer. Therefore, we have
\BQn
\sum_{\substack{\alpha \in \bN^n \\ |\alpha|=m }} \alpha !
\leq m! \sum_{\substack{\alpha \in \bN^n \\ |\alpha|=m }} 1. 
\EQn

Secondly, note that $\sum_{\substack{\alpha \in \bN^n \\ |\alpha|=m }} 1$ 
is just the number of distinct $\alpha\in \bN^n$ with $|\alpha|=m$, 
which is the same as the number of distinct monomials in $n$ free commutative 
variables of degree $m$. Since the latter is well-known to be the binomial 
$\binom{m+n-1}{m}$, we have 
\BQn
\sum_{\substack{\alpha \in \bN^n \\ |\alpha|=m }} \alpha !
\leq m! \binom{m+n-1}{m} 
=\frac{(m+n-1)!}{(n-1)!}. 
\EQn
\epfv

\underline{\it Proof of Proposition \ref{cvg-MainThm-2}}:
By Eq.\,$(\ref{E-Q-1})$ , we know that, 
\begin{align}\label{cvg-E-Q-1}
Q(z)=\sum_{m\geq 1} \frac {\Delta^m P^{m+1}} {2^m m!(m+1)!}.
\end{align}

To show the proposition, it will be enough to show 
the infinite series above converges absolutely
over $B(0, r)$ for any $r<r_0$. 

We first give an upper bound for the general terms 
in the series Eq.\,$(\ref{cvg-E-Q-1})$ over $B(0, r)$.

Consider
\begin{align}
\Delta^m P^{m+1}&= ( \sum_{i=1}^n D_i^2 )^m P^{m+1} 
=\sum_{\substack{\alpha\in \bN^n \\ |\alpha|=m}}
\frac{m!}{\alpha !} D^{2\alpha} P^{m+1}. 
\end{align}

Therefore, we have
\allowdisplaybreaks{
\begin{align*}
|\Delta^m P^{m+1}(a)|
& \leq \sum_{\substack{\alpha\in \bN^n \\ |\alpha|=m }}
\frac{m!}{\alpha !} |D^{2\alpha} P^{m+1}(a)| \\
\intertext{Applying Lemma \ref{cvg-L1} with $\alpha$ replaced by $2\alpha$:}
& \leq \sum_{\substack{\alpha \in \bN^n \\ |\alpha|=m }}
\frac{m!}{\alpha !} \frac{(2\alpha)!}{r^{2m}} (2r)^{d(m+1)} |P|^{m+1} \\
\intertext{Noting that $(2\alpha)! \leq [(2\alpha)!!]^2
=2^{2m}(\alpha!)^2$:}
& \leq \sum_{\substack{\alpha \in \bN^n \\ |\alpha|=m }}
\frac{m!}{\alpha !} \frac{2^{2m}(\alpha !)^2}{ r^{2m} } 
(2r)^{d(m+1)} |P|^{m+1} \\
&= m! 2^{2m+d(m+1)} r^{d(m+1)-2m} |P|^{m+1} \sum_{\substack{\alpha \in \bN^n \\ |\alpha|=m }} \alpha ! \\
\intertext{Applying Lemma \ref{cvg-L2}:} 
&=\frac{m! (m+n-1)! 2^{2m+d(m+1)} r^{d(m+1)-2m} |P|^{m+1}}{(n-1)!}.
\end{align*}
Therefore, for any $m\geq 1$, we have
\begin{align}\label{cvg-Am}
\left | \frac{\Delta^m P^{m+1}} {2^m m!(m+1)!}  \right |
 \leq \frac{2^{m+d(m+1)} r^{d(m+1)-2m} |P|^{m+1} (m+n-1)! }{(m+1)! (n-1)!}.
\end{align}
}
For any $m\geq 1$, let $A_m$ be the right hand side of Eq.\,$(\ref{cvg-Am})$  above. 
Then, by a straightforward calculation, we see that the ratio
\begin{align}
\frac{A_{m+1}}{A_m}= \frac{m+n}{m+2}\, 2^{d+1} r^{d-2}|P|.
\end{align}

Since $r<r_0=(2^{d+1}|P|)^{\frac 1{2-d}}$, it is easy to see that
\begin{align*}
\lim_{m\to +\infty} \frac{A_{m+1}}{A_m} = 2^{d+1} r^{d-2}|P|<1. 
\end{align*}

Therefore, by the comparison test, the infinite series 
in Eq.\,$(\ref{cvg-E-Q-1})$ 
converges absolutely and uniformly 
over the open ball $B(0, r)$. 
\epfv

\renewcommand{\theequation}{\thesection.\arabic{equation}}
\renewcommand{\therema}{\thesection.\arabic{rema}}
\setcounter{equation}{0}
\setcounter{rema}{0}

\section{\bf Self-Inverting Formal Power Series}\label{slf}

Note that, by the definition of inversion pairs 
(see page \pageref{Inv-pair}), 
 $Q\in \bCzz$ is the inversion pair 
of $P\in \bCzz$ iff $P$ is the inversion pair of $Q$. 
In other words, the relation that $Q$ and $P$ are 
inversion pair of each other in some sense 
is a duality relation. Naturally, one may ask, 
for which $P(z)$, it is self-dual 
or self-inverting? In this section, we 
discuss this special family of 
polynomials or formal power series. 

Another purpose of this section 
is to draw the reader's attention to
the problem of classification of (HN) 
self-inverting polynomials (see Open Problem \ref{slf-prob}). 
Even though the classification of HN polynomials 
seems to be out of reach at the current time, 
we believe that the classification of (HN) 
self-inverting polynomials is much 
more approachable.

\begin{defi}\label{slf-Def}
A formal power series $P(z)\in \bCzz$ 
with $o(P(z))\ge 2$ 
and $(\Hes P)(0)$ nilpotent is said to be 
{\it self-inverting} if its inversion pair 
$Q(z)=P(z)$.
\end{defi}

Following the terminology introduced 
in \cite{B2}, we say a formal map $F(z)=z-H(z)$ 
with $H(z)\in \bCzz^{\times n}$ and $o(H(z))\geq 1$ 
is a {\it quasi-translation} if $j(F)(0) \ne 0$ and
its formal inverse map is given 
by $G(z)=z+H(z)$. 

Therefore, for any $P(z)\in \bC[[z]]$ with 
$o(P(z))\geq 2$ and $(\Hes P)(0)$ nilpotent, 
it is self-inverting iff the associated symmetric 
formal map $F(z)=z-\nabla P(z)$ 
is a quasi-translation. 

For quasi-translations, the following general 
result has been proved in Proposition $1.1$ of  
\cite{B2} for polynomial quasi-translations. 

\begin{propo} \label{B2-P}
A formal map $F(z)=z-H(z)$ with $o(H)\geq 1$ 
and $JH(0)$ nilpotent is a quasi-translation
if and only if $JH \cdot H =0$.
\end{propo}

Even though the proposition above was proved in 
\cite{B2} only in the setting of polynomial 
maps, the proof given there works 
equally well for formal quasi-translations 
under the condition that $JH(0)$ is nilpotent. 
Since it has also been shown in Proposition $1.1$ 
in \cite{B2} that, for any polynomial 
quasi-translations $F(z)=z-H(z)$, 
$JH(z)$ is always nilpotent,  
so the condition that $JH(0)$ is 
nilpotent in the proposition above 
does not put any extra restriction 
for the case of polynomial 
quasi-translations.  

From Proposition \ref{B2-P} above, we immediately have 
the following criterion for 
self-inverting formal power series.

\begin{propo} \label{slf-MainThm-1}
For any $P(z)\in \bCzz$ with $o(P)\geq 2$ and 
$(\Hes P)(0)$ nilpotent, it is self-inverting 
if and only if $\la\nabla P, \nabla P\ra =0$.
\end{propo}

\pf Since $o(P)\geq 2$ and 
$(\Hes P)(0)$ is nilpotent, by Proposition \ref{B2-P}, 
we see that, $P(z)\in \bCzz$ is self-inverting iff 
$J(\nabla P)\cdot \nabla P=(\Hes P) \cdot \nabla P=0$.
But, on the other hand,  it is easy to check that, 
for any $P(z)\in \bCzz$, we have the following identity:
\begin{align*}
(\Hes P) \cdot \nabla P= \frac 12 \nabla \la \nabla P, \nabla P\ra.
\end{align*}
Therefore, $(\Hes P) \cdot \nabla P=0$ 
iff $\nabla \la \nabla P, \nabla P \ra=0$, and iff
$\la \nabla P, \nabla P\ra=0$ because 
$o(\la \nabla P, \nabla P\ra)\geq 2$.
\epfv

\begin{corol}\label{slf-CC2}
For any $P(z)\in \bCzz$ with $o(P)\geq 2$ and 
$(\Hes P)(0)$ nilpotent, if it is self-inverting, 
then so is $P^m(z)$ for any $m\geq 1$.
\end{corol}
\pf Note that, for any $m\ge 2$, we have $o(P^m(z))\ge 2m >2$ and 
$(\Hes P)(0)=0$. Then, the corollary follows 
immediately from Proposition \ref{slf-MainThm-1} 
and the following general identity:
\begin{align}
\la \nabla P^m, \nabla P^m \ra= m^2P^{2m-2} 
\la \nabla P, \nabla P\ra.
\end{align}
\epfv

\begin{corol}\label{slf-C2}
For any harmonic formal power series $P(z)\in \bCzz$ with $o(P)\geq 2$ and 
$(\Hes P)(0)$ nilpotent, it is self-inverting 
iff $\Delta P^2=0$.
\end{corol}
\pf This follows immediately from Proposition 
\ref{slf-MainThm-1} and the following  
general identity:
\begin{align}
\Delta P^2=2(\Delta P)P + 2\la \nabla P, \nabla P\ra.
\end{align}
\epfv

\begin{propo}\label{slf-MainThm-2}
Let $P(z)$ be a harmonic self-inverting 
formal power series. 
Then, for any $m\geq 1$, $P^m$ is HN. 
\end{propo}

\pf First, we use the mathematical induction on 
$m\geq 1$ to show that $\Delta P^m=0$ for any $m\geq 1$. 

The case of $m=1$ is given. For any $m\geq 2$, consider
\allowdisplaybreaks{
\begin{align*}
\Delta P^m &= \Delta (P\cdot P^{m-1})\\
&=(\Delta P) P^{m-1}+ P(\Delta P^{m-1})
+ 2\la \nabla P, \nabla P^{m-1}\ra \\
&=(\Delta P) P^{m-1}+ P(\Delta P^{m-1})
+2(m-1) P^{m-2} \la \nabla P, \nabla P\ra.
\end{align*} }
Then, by the mathematical induction assumption and Proposition  
\ref{slf-MainThm-1}, we get $\Delta P^m=0$.

Secondly, for any fixed $m\geq 1$ and $d\geq 1$, we have
\begin{align*}
\Delta^d [(P^m)^d] &= \Delta^{d-1} (\Delta P^{dm})=0.
\end{align*}
Then, by the criterion in Proposition \ref{Crit-1}, 
$P^m$ is HN.
\epfv

\begin{exam}
Note that, in Section $5.2$ of \cite{HNP}, 
a family of self-inverting HN formal power series 
has been constructed as follows.

Let $\Xi$ be any non-empty subset of $\bC^n$ 
such that, for any $\alpha, \beta\in \Xi$, 
$\la \alpha, \beta \ra=0$. 
Let $\cA$ be the completion 
of the subalgebra of $\bCzz$ generated by 
$h_\alpha (z)\!:=\la \alpha, z \ra$ 
$(\alpha\in \Xi)$, i.e. 
$\cA$ is the set of all formal power series in 
$h_\alpha (z)$ $(\alpha\in \Xi)$ over $\bC$.
Then it is straightforward 
to check $($or see Section $5.2$ of \cite{HNP} for details$)$
that any element $P(z)\in \cA$ 
is HN and self-inverting. 
\end{exam}

It is unknown if all HN self-inverting polynomials 
or formal power series can be obtained 
by the construction above. More generally, 
we believe the following open problem 
is worth investigating. 

\begin{prob}\label{slf-prob}
$(a)$ Decide whether or not all self-inverting polynomials 
or formal power series are HN.

$(b)$ Classify all $($HN$)$ self-inverting  
polynomials and formal power series.
\end{prob}

Finally, let us point out that, for any 
self-inverting $P(z)\in \bCzz$, the 
deformed inversion pair $Q_t(z)$ (not just $Q(z)=Q_{t=1}(z)$) 
is also same as $P(z)$. 

\begin{propo}\label{Qt=P} 
Let $P(z)\in \bCzz$ with $o(P)\geq 2$ and 
$(\Hes P)(0)$ nilpotent. Then $P(z)$ 
is self-inverting 
if and only if $Q_t(z)=P(z)$. 
\end{propo}

\pf First, let us point out the following observations.

Let $t$ be a formal central parameter 
and $F_t(z)=z-t \nabla P(z)$ as before. 
Since $o(P)\geq 2$ and $(\Hes P)(0)$ is nilpotent, 
we have $j(F_t)(0)=1$. Therefore, $F_t(z)$ 
is an automorphism of the 
algebra $\bC[t][[z]]$ of 
formal power series of $z$ 
over $\bC[t]$. Since the inverse map of $F_t(z)$ 
is given by $G_t(z)=z+t\nabla Q_t(z)$, we see 
that $Q_t(z)\in \bC[t][[z]]$. Therefore, 
for any $t_0\in \bC$, $Q_{t=t_0}(z)$ 
makes sense and lies in $\bC[[z]]$. 
Furthermore, by the uniqueness of inverse 
maps, it is easy to see that the 
inverse map of $F_{t_0}=z-t_0 \nabla P$ 
of $\bC[t][[z]]$ is given by 
$G_{t_0}(z)=z+ t_0 \nabla Q_{t=t_0}$. 
Therefore the inversion pair of 
$t_0 P(z)$ is given by 
$t_0  Q_{t=t_0}(z)$.

With the notation and observations above, 
by choosing $t_0=1$, we have $Q_{t=1}(z)=Q(z)$ 
and the $(\Leftarrow)$ part of 
the proposition follows immediately. 
Conversely, for any $t_0\in \bC$, we have 
$\la \nabla (t_0 P), \nabla (t_0 P)\ra
=t_0^2 \la \nabla P, \nabla P\ra$. Then, 
by Proposition \ref{slf-MainThm-1},  
$t_0 P(z)$ is self-inverting and its inversion pair 
$t_0 Q_{t=t_0}(z)$ is same as $t_0 P(z)$, i.e.  
$t_0 Q_{t=t_0}(z)=t_0 P(z)$. 
Therefore, we have $ Q_{t=t_0}(z)=P(z)$ 
for any $t_0\in \bC^\times$. 
But on the other hand, we have 
$Q_t(z)\in \bC[t][[z]]$ as pointed 
above, i.e. the coefficients of all 
monomials of $z$ in $Q_t(z)$ are 
polynomials of $t$, hence we must have 
$Q_t(z)=P(z)$ which is the $(\Rightarrow)$ part of 
the proposition. 
\epfv

\renewcommand{\theequation}{\thesection.\arabic{equation}}
\renewcommand{\therema}{\thesection.\arabic{rema}}
\setcounter{equation}{0}
\setcounter{rema}{0}

\section{\bf The Vanishing Conjecture over 
Fields of Positive Characteristic} \label{pVC}

It is well-known that the JC may fail when $F(z)$ is not 
a polynomial map  (e.g. $F_1(z_1, z_2)=e^{-z_1};$ 
$F_2(z_1, z_2)=z_2 e^{z_1}$). It also fails
badly over fields of positive characteristic 
even in one variable case 
(e.g. $F(x)=x-x^p$ over a field of characteristic $p>0$). 
However, the situation for the VC over fields 
of positive characteristic is dramatically 
different from the JC even through these two conjectures 
are equivalent to each other  over fields of 
characteristic zero. Actually, as we will show 
in the proposition below, the VC over fields 
of positive characteristic holds for 
any polynomials (not even necessarily HN) and also 
for any HN formal power series.

\begin{propo} \label{pVC-MainThm}
Let $k$ be a field of characteristic  $p>0$. Then

$(a)$ For any polynomial $P(z)\in k[z]$ $($not necessarily 
homogeneous  nor HN$)$ of degree $d\geq 1$, \, 
$\Delta^m P^{m+1} =0$ for any $m \geq \frac {d(p-1)}{2}$.

$(b)$ For any HN formal power series $P(z)\in k[[z]]$, i.e.
$\Delta^m P^{m}=0$ for any $m\geq 1$, we have,   
$\Delta^m P^{m+1}=0$ for any $m\geq p-1$.
\end{propo}

In other words, over the fields of positive 
characteristic, the VC holds even for HN 
formal power series $P(z)\in k[[z]]$; 
while for polynomials, it holds even 
without the HN condition nor 
any other conditions. 

\pf The main reason that 
the proposition above holds is because of
the following simple fact due to the Leibniz rule 
and positiveness of the characteristics of the base field 
$k$, namely, for $m\geq 1$,  $u(z), v(z)\in k[[z]]$ and 
any differential operator $\Lambda$ of $k[z]$, 
we have
\begin{align}\label{pVC-MainThm-pe1}
\Lambda (u^{mp} v)=u^{mp}\Lambda v.
\end{align}
 
Now let $P(z)$ be any polynomial or formal series 
as in the proposition. For any $m\geq 1$, 
write $m+1=q_m p+r_m$ with $q_m, r_m\in \bZ$ and 
$0\leq r_m \leq p-1$. 
Then by Eq.\,$(\ref{pVC-MainThm-pe1})$ , we have 
\begin{align}\label{pVC-MainThm-pe2}
\Delta^m P^{m+1}= \Delta^m (P^{q_mp}P^{r_m})=P^{q_mp} \Delta^m P^{r_m}.
\end{align}

If $P(z)$ is a polynomial of degree $d \geq  1$, 
we have $\Delta^m P^{r_m} =0$ when  
$m \geq \frac{d(p-1)}2 $, since in this case 
$2m>\deg (P^{r_m})$. If $P(z)$ is a HN 
formal power series, we have 
$\Delta^m P^{r_m}=0$ when $m \geq p-1 \ge r_m$. 
Therefore, $(a)$ and $(b)$ in the proposition 
follow from Eq.\,$(\ref{pVC-MainThm-pe2})$  
and the observations above.
\epfv

One interesting question is whether or not
the VC fails (as the JC does) for any HN formal 
power series $P(z)\in \bC[[z]]$ but 
$P(z)\not \in \bC[z]$?  To our best knowledge, 
no such counterexample has been known yet. 
We here put it as an open problem.

\begin{prob}\label{pVC-prob}
Find a HN formal power series  $P(z)\in \bC[[z]]$ but 
$P(z)\not \in \bC[z]$,  if there are any, such that  
the VC fails for $P(z)$. 
\end{prob}

One final remark about Proposition \ref{pVC-MainThm}
is as follows. Note that the crucial fact used 
in the proof is that any differential operator $\Lambda$ 
of $k[z]$ commutes with the multiplication operator 
by the $p^{th}$ power of any element of $k[[z]]$. 
Then, by a parallel argument as in the proof of 
Proposition \ref{pVC-MainThm}, it is easy to see that
the following more general result also holds.

\begin{propo}
Let $k$ be a field of characteristics $p>0$ and 
$\Lambda$ a differential operator of $k[z]$. 
Let $f\in k[[z]]$. Assume that, 
for any $1\leq m\leq p-1$, 
there exists $N_m>0$ such that 
$\Lambda^{N_m} f^m=0$.
Then, we have $\Lambda^{m} f^{m+1}=0$ 
when $m>>0$.

In particular, if $\Lambda$ strictly 
decreases the degree of polynomials. 
Then, for any polynomial $f\in k[z]$, 
we have $\Lambda^{m} f^{m+1}=0$  
when $m>>0$.
\end{propo}

\renewcommand{\theequation}{\thesection.\arabic{equation}}
\renewcommand{\therema}{\thesection.\arabic{rema}}
\setcounter{equation}{0}
\setcounter{rema}{0}

\section{\bf A Criterion of Hessian Nilpotency for Homogeneous Polynomials}
\label{crn}

Recall that $\la \cdot, \cdot \ra$ denotes the standard $\bC$ 
bilinear form of $\bC^n$. For any $\beta\in \bC^n$, 
we set $h_\beta(z)\!:=\la \beta, z\ra$ and 
$\beta_D\!:=\la \beta, D\ra$. 

The main result of this section 
is the following criterion of 
Hessian nilpotency for homogeneous polynomials. 
Considering the criterion given in 
Proposition \ref{Crit-1}, it is somewhat 
surprising but the proof turns out to be very simple.   

\begin{theo}\label{crn-MainThm}
For any $\beta \in \bC^n$ and homogeneous  polynomial $P(z)$ 
of degree $d\geq 2$, set $P_\beta (z)\!:= \beta_D^{d-2} P(z)$. 
Then, we have 
\begin{align}\label{crn-MainThm-e1}
{\Hes P_\beta}=(d-2)! \, (\Hes P)(\beta).
\end{align}
In particular, $P(z)$ is HN 
iff, for any $\beta\in \bC^n$, $P_\beta (z)$ is HN.
\end{theo}

To prove the theorem, we need first the following lemma.

\begin{lemma}\label{crn-L1}
Let $\beta\in\bC^n$ and $P(z)\in \bCz$ homogeneous  
of degree $N \geq 1$. Then
\begin{align}\label{crn-L1-e1}
\beta_D^N P(z) = N!\,P(\beta).
\end{align}
\end{lemma}
\pf Since both sides of Eq.\,$(\ref{crn-L1-e1})$  
are linear on $P(z)$, we may assume 
$P(z)$ is a monomial, say $P(z)=z^{\bf a}$ 
for some ${\bf a}\in \bN^n$ 
with $|{\bf a}|=N$. 

Consider
\allowdisplaybreaks{
\begin{align*}
\beta_D^N P(z) &=(\sum_{i=1}^n \beta_iD_i )^N z^{\bf a}
=\sum_{\substack{ {\bf k}\in \bN^n \\ |{\bf k}|=N} }^n 
\frac {N!}{{\bf k}!} \beta^{\bf k}D^{\bf k} z^{\bf a} \\ 
&=\frac {N!}{{\bf a}!} \beta^{\bf a} D^{\bf a} z^{\bf a}
=N! \beta^{\bf a} 
= N!P(\beta).
\end{align*}
}
\epfv

\underline{\it Proof of Theorem \ref{crn-MainThm}}: 
We consider 
\allowdisplaybreaks{
\begin{align*}
\Hes P_\beta (z) &= 
\left ( \frac{\p^2 (\beta_D^{d-2} P) }{\p z_i\p z_j}(z) \right )_{n\times n}  
=\left ( \beta_D^{d-2} \frac{\p^2 P}{\p z_i\p z_j}(z) \right )_{n\times n}  \\
\intertext{Applying Lemma \ref{crn-L1} to $\frac{\p^2 P}{\p z_i\p z_j}(z)$:}
&=(d-2)! \left ( \frac{\p^2 P}{\p z_i\p z_j}(\beta) \right )_{n\times n}  
=(d-2)! \, (\Hes P)(\beta).
\end{align*} }
\epfv

Let $\{e_i\,|\, 1\leq i\leq n\}$ be the standard basis of $\bC^n$. 
Applying the theorem above to $\beta=e_i$ $(1\leq i\leq n)$, 
we have the following corollary, which was first 
proved by M. Kumar \cite{Ku}.

\begin{corol}\label{crn-Corol}
For any homogeneous HN polynomial $P(z)\in \bC[z]$ 
of degree $d\ge 2$, $D_i^{d-2}P(z)$  
$(1\leq i\leq n)$ are also HN. 
\end{corol}

The reason that we think the criteria given in 
Theorem \ref{crn-MainThm} and Corollary \ref{crn-Corol}
interesting is that, $P_\beta(z)=\beta_D^{d-2}P(z)$ 
is homogeneous of degree $2$, and it is much easier 
to decide whether a homogeneous polynomial 
of degree $2$ is HN or not. More precisely, 
for any homogeneous polynomial $U(z)$  
of degree $2$, there exists a unique symmetric 
$n\times n$ matrix $A$ such that 
$U(z)=z^\tau A z$. Then it is easy to 
check that $\Hes U(z)=2A.$ Therefore, 
$U(z)$ is HN iff the symmetric matrix $A$
is nilpotent.

Finally we end this section with the following 
open question on the criterion given 
in Proposition \ref{Crit-1}.

Recall that Proposition \ref{Crit-1} 
was proved in \cite{HNP}. We now sketch the argument.

For any $m\geq 1$, we set 
\begin{align}
u_m(P)&=\text{Tr\,}  \text{Hes}^m(P), \label{Def-u} \\
v_m(P)&=\Delta^m P^m.\label{Def-v}
\end{align}

For any $k\geq 1$, we define ${\mathcal U}_k(P)$ 
(resp.\,${\mathcal V}_k(P)$) to
be the ideal in $\bC[[z]]$ generated by  
$\{u_m(P) |1\leq m\leq k\}$ 
(resp.\,$\{ v_m(P) | 1\leq m\leq k \}$)
and all their partial derivatives 
of any order. Then it has been shown 
(in a more general setting) 
in Section $4$ in \cite{HNP} that
$\mathcal U_k(P)=\mathcal V_k(P)$
for any $k\geq 1$.   

It is well-known 
in linear algebra that, if $u_m(P(z))=0$ 
when $m>>0$, then $\Hes P$ is 
nilpotent and $u_m(P)=0$ 
for any $m\geq 1$. One natural 
question is whether or not this 
is also the case for 
the sequence $\{v_m(P)\,|\, m\geq 1\}$.
More precisely, we believe the following 
conjecture which was
proposed in \cite{HNP} is worth investigating.

\begin{conj}\label{Conj-4.3-HNP}
Let $P(z)\in \bC [[z]]$ with $o(P(z))\geq 2$. 
If $\Delta^m P^m(z)=0$ for $m>>0$, then $P(z)$ is HN.
\end{conj}

\renewcommand{\theequation}{\thesection.\arabic{equation}}
\renewcommand{\therema}{\thesection.\arabic{rema}}
\setcounter{equation}{0}
\setcounter{rema}{0}

\section{\bf Some Results on Symmetric Polynomial Maps}
\label{symmc}

Let $P(z)$ be any formal power series with $o(P(z))\geq 2$ 
and $(\Hes P)(0)$ nilpotent, and $F(z)$ and $G(z)$ as before. 
Set
\begin{align}
\sigma_2\!:&=\sum_{i=1}^n z_i^2,\\ 
f(z)\!:&=\frac 12 \sigma_2- P(z).
\end{align} 
Professors Mohan Kumar \cite{Ku} and 
David Wright \cite{Wr3} once asked 
how to write $P(z)$ and $f(z)$ 
in terms of $F(z)$? More precisely, 
find $U(z), V(z) \in \bCzz$ 
such that 
\begin{align}
U(F(z))&=P(z),\label{UFP} \\
V(F(z))&=f(z).\label{VFP}
\end{align}

In this section, we first derive 
in Proposition \ref{symmc-MainThm-1}
some explicit formulas for $U(z)$ and $V(z)$, and 
also for $W(z)\in \bCzz$ such that 
\begin{align}\label{WFP}
W(F(z))=\sigma_2(z).
\end{align}
We then show in Theorem \ref{symmc-MainThm-2} that, 
when $P(z)$ is a HN polynomial, the VC holds for $P$ or 
equivalently, the JC holds for the associated 
symmetric polynomial map $F(z)=z-\nabla P$, 
iff one of $U$, $V$ and $W$ 
is polynomial.

Let $t$ be a central parameter and $F_t(z)=z-t\nabla P$. 
Let $G_t(z)=z+t\nabla Q_t$ be the formal 
inverse of $F_t(z)$ as before. 
We set 
\allowdisplaybreaks{
\begin{align}
f_t(z)\!:&= \frac 12 \sigma_2- tP(z),\\
U_t(z)\!:&=P( G_t(z) ), \label{Def-Ut}  \\
V_t(z)\!:&= f_t(G_t(z)), \label{Def-Vt}\\
W_t(z)\!:&=\sigma_2(G_t(z)). \label{Def-Wt}
\end{align} }

Note first that, under the conditions that 
$o(P(z)) \ge 2$ and $(\Hes P)(0)$ 
is nilpotent, we have 
$G_t(z)\in \bC[t][[z]]^{\times n}$ 
as mentioned in the proof of Proposition 
\ref{Qt=P}. Therefore, we have 
$U_t(z), V_t(z), W_t(z) \in \bC[t][[z]]$, 
and $U_{t=1}(z)$, $V_{t=1}(z)$ and $W_{t=1}(z)$ 
all make sense. Secondly, from the 
definitions above, we have
\begin{align}
W_t(z)&=2V_t(z)+2tU_t(z),\label{W=VU} \\
F_t(z)&=\nabla f_t(z),\\
f_{t=1}(z)&= f(z).
\end{align}

\begin{lemma}\label{lemma-0}
With the notations above, we have
\begin{align}
P(z)&=U_{t=1}(F(z)),\label{lemma-0-e1} \\
f(z)&=V_{t=1}(F(z)),\label{lemma-0-e2}\\
\sigma_2(z)&=W_{t=1}(F(z)). \label{lemma-0-e3}
\end{align}
In particular, $f(z)$, $P(z)$ and $\sigma_2(z)$ lie in 
$\bC[F]$ iff $U_{t=1}(z)$, $V_{t=1}(z)$ and
$W_{t=1}(z)$ lie in $\bC[z]$.
\end{lemma}

In other words, by setting $t=1$, $U_t$, $V_t$ and $W_t$ will 
give us $U$, $V$ and $W$ in Eqs.\,$(\ref{UFP})$--$(\ref{WFP})$, 
respectively.

\pf From the definitions of $U_t(z)$, $V_t(z)$ and $W_t(z)$
 (see Eqs.\,$(\ref{Def-Ut})$--$(\ref{Def-Wt})$, we have
\begin{align*}
P(z)&=U_{t}(F_t(z)),\\
f_t(z)&=V_{t}(F_t(z)),\\
\sigma_2(z)&=W_{t}(F_t(z)). 
\end{align*}
By setting $t=1$ in the equations above 
and noticing that $F_{t=1}(z)=F(z)$, we get 
Eqs.\,$(\ref{lemma-0-e1})$--$(\ref{lemma-0-e3})$. 
\epfv

For $U_t(z)$, $V_t(z)$ and $W_t(z)$, 
we have the following explicit formulas in 
terms of the deformed inversion pair $Q_t$ of $P$.

\begin{propo}\label{symmc-MainThm-1}
For any formal power series $P(z)\in \bCzz$ $($not necessarily HN$)$ 
with $o(P(z))\geq 2$ and $(\Hes P)(0)$ nilpotent, we have
\allowdisplaybreaks{
\begin{align}
U_t(z) &=Q_t+ t \, \frac{\p Q_t}{\p t},
 \label{U-t} \\
V_t(z) &=\frac 12\sigma_2 + t (z\frac {\p Q_t}{\p z} -Q_t),  \label{V-t}\\
W_t(z) &=\sigma_2 + 2t z\frac {\p Q_t}{\p z} +2t^2
\frac{\p Q_t}{\p t}. \label{W-t}
\end{align} }
\end{propo}

\pf Note first that, Eq.\,(\ref{W-t}) follows directly from 
Eqs.\,$(\ref{U-t})$, $(\ref{V-t})$ and $(\ref{W=VU})$.

To show Eq.\,(\ref{U-t}),
by Eqs.\,$(3.4)$ and $(3.6)$ in \cite{Burgers}, 
we have 
\begin{align}\label{E3.4-Z1}
U_t(z)=P(G_t)
 =Q_t+\frac t2 \la \nabla Q_t, \nabla Q_t \ra 
 =Q_t+ t\, \frac{\p Q_t}{\p t}. 
\end{align}

To show Eq.\,(\ref{V-t}),  we consider
\allowdisplaybreaks{
\begin{align*}
V_t(z)&=f_t(G_t)\\
&=\frac 12 \la z+t\nabla Q_t(z), z+t\nabla Q_t(z)\ra -tP(G_t)\\
&=\frac 12\sigma_2 +t \la z, \nabla Q_t(z)\ra +\frac {t^2}2 
\la \nabla Q_t, \nabla Q_t \ra 
-tP(G_t) \\
\intertext{By Eq.\,(\ref{E3.4-Z1}), 
substituting $ Q_t+\frac t2 \la \nabla Q_t, \nabla Q_t \ra$
for $P(G_t)$:}
&=\frac 12\sigma_2 +t \la z, \nabla Q_t(z)\ra - t Q_t(z)\\ 
&=\frac 12\sigma_2 + t (z\frac {\p Q_t}{\p z} -Q_t). 
\end{align*} }
\epfv

When $P(z)$ is homogeneous and HN, we have the following more 
explicit formulas which in particular give solutions 
to the questions raised by Professors Mohan Kumar 
and David Wright.

\begin{corol}\label{symmc-MainCor-1}
For any homogeneous HN polynomial $P(z)$ of degree $d\geq 2$, we have
\allowdisplaybreaks{
\begin{align}
U_t(z) & =\sum_{m=0}^\infty \frac {t^m}{2^m (m!)^2} \Delta^m P^{m+1}(z) 
 \label{U-t-2} \\
V_t(z) &= \frac 12 \sigma_2+ 
\sum_{m=0}^\infty \frac{(d_m-1)t^{m+1}}{2^m m!(m+1)!}\, \Delta^m P^{m+1}(z)\, , 
\label{V-t-2} \\
W_t(z) &= \sigma_2+
\sum_{m=0}^\infty \frac {(d_m+m)t^{m+1}}{2^{m-1} m! (m+1)!}\, \Delta^m P^{m+1}(z) \, ,\label{W-t-2}
\end{align} }
where $d_m=\deg\, (\Delta^m P^{m+1})=d(m+1)-2m$ $(m\geq 0)$. 
\end{corol}

\pf We give a proof for Eq.\,(\ref{U-t-2}). 
Eqs.\,$(\ref{V-t-2})$ can be proved similarly. 
$(\ref{W-t-2})$ follows directly from  
Eqs.\,Eq.\,(\ref{U-t-2}), \,$(\ref{V-t-2})$  
and $(\ref{W=VU})$.

By combining Eqs.\,(\ref{U-t}) and (\ref{E-Q-1}), we have
\allowdisplaybreaks{
\begin{align*}
U_t(z)&= \sum_{m=0}^\infty \frac {t^{m} \Delta^m P^{m+1}(z)}{2^m m!(m+1)!}  +
 \sum_{m=1}^\infty \frac { m t^m  \Delta^m P^{m+1}(z) }{2^m m!(m+1)!}  \\
&= P(z)+ \sum_{m=1}^\infty \frac {t^m}{2^m (m!)^2} \Delta^m P^{m+1}(z) \\
&=\sum_{m=0}^\infty \frac {t^m}{2^m (m!)^2} \Delta^m P^{m+1}(z). 
\end{align*} }
Hence, we get Eq.\,(\ref{U-t-2}). 
\epfv

One consequence of the proposition 
above is the following result 
on symmetric polynomials maps.

\begin{theo}\label{symmc-MainThm-2}
For any HN polynomial $P(z)$ $($not necessarily homogeneous$)$ 
with $o(P)\geq 2$,  
the following statements are equivalent: 
\begin{enumerate}
\item The VC holds for $P(z)$.
\item $P(z)\in \bC[F]$.
\item $f(z)\in \bC[F]$.
\item $\sigma_2(z)\in \bC[F]$.
\end{enumerate}
\end{theo}

Note that, the equivalence of 
the statements $(1)$ and $(3)$
was first proved by Mohan Kumar 
(\cite{Ku}) by a different method. \\

\pf Note first that, by Lemma \ref{lemma-0}, 
it will be enough to show that, 
$\Delta^m P^{m+1}=0$ when $m>>0$ iff 
one of $U_t(z)$, $V_t(z)$ and $W_t(z)$ is a
polynomial in $t$ with coefficients in 
$\bCz$. Secondly, when $P(z)$ is 
homogeneous, the statement above 
follows directly from 
Eqs.\,$(\ref{U-t-2})$--$(\ref{W-t-2})$.

To show the general case, 
for any $m\ge 0$ and $M_t(z) \in \bC[t][[z]]$, 
we denote by $[t^m]( M_t(z) )$ the coefficient
of $t^m$ when we write $M_t(z)$ 
as a formal power series of $t$ 
with coefficients in $\bCzz$.
Then, from  Eqs.\,$(\ref{U-t})$--$(\ref{W-t})$ 
and Eq.\,$(\ref{E-Q-1})$, it is straightforward  
to check that the coefficients of $t^m$ $(m\geq 1)$ 
in $U_t(z)$, $V_t(z)$ and $W_t(z)$ 
are given as follows.
\allowdisplaybreaks{
\begin{align}
[t^m]( U_t(z) ) = \frac {\Delta ^m P^{m+1}}{2^m(m!)^2}, \label{Cof-Ut} 
\end{align}
\begin{align}
[t^m]( V_t(z) ) = \frac {1}{2^{m-1}(m-1)! m!}
\left (z \frac{\p}{\p z} (\Delta^{m-1} P^m )-\Delta^{m-1} P^m  \right), 
\label{Cof-Vt} 
\end{align}
\begin{align}
[t^m]( W_t(z) ) = \frac {1}{2^{m-2}(m-1)! m!}
\left (z \frac{\p}{\p z} (\Delta^{m-1} P^m )
+ (m-1)\Delta^{m-1} P^m  \right).\label{Cof-Wt}
\end{align} }

From Eq.\,$(\ref{Cof-Ut})$, we immediately have 
$(1)\Leftrightarrow (2)$. 
To show the equivalences 
$(1)\Leftrightarrow (3)$ 
and $(1)\Leftrightarrow (4)$, 
note first that $o(P)\geq 2$, 
so $o(\Delta^{m-1} P^m) \geq 2$ 
for any $m\geq 1$. While, on the other hand, 
for any polynomial $h(z)\in \bC[z]$ with 
$o(h(z))\geq 2$, we have, 
$h(z)=0$ iff $( z\frac{\p}{\p z}-1 )h(z)=0$, 
and iff $( z\frac{\p}{\p z}+ (m-1) )h(z)=0$ 
for some $m\geq 1$. This is simply because that, 
for any monomial $z^\alpha$ $(\alpha\in \bN^n)$,  
we have $( z\frac{\p}{\p z}-1 )z^\alpha =(|\alpha|-1)z^\alpha$ 
and $( z\frac{\p}{\p z}+ (m-1) )z^{\alpha}
=(|\alpha|+(m-1))z^\alpha$. From this general fact, we see that 
$(1)\Leftrightarrow (3)$ follows from Eq.\,$(\ref{Cof-Vt})$
and $(1)\Leftrightarrow (4)$ from Eq.\,$(\ref{Cof-Wt})$.
\epfv

\renewcommand{\theequation}{\thesection.\arabic{equation}}
\renewcommand{\therema}{\thesection.\arabic{rema}}
\setcounter{equation}{0}
\setcounter{rema}{0}

\section{\bf A Graph Associated with Homogeneous HN Polynomials}
\label{grf}

In this section, we would like to draw the reader's 
attention to a graph $\cG(P)$ assigned to 
each homogeneous harmonic polynomials $P(z)$. 
The graph $\cG(P)$ was first proposed by the author 
and later was further studied by R. Willems  
in his master thesis \cite{Roel} 
under direction of Professor A. van den Essen. 
The introduction of the graph $\cG(P)$ is mainly 
motivated by a criterion of Hessian nilpotency 
given in \cite{HNP} 
(see also Theorem \ref{Crit-2} below), 
via which one hopes more 
necessary or sufficient conditions 
for a homogeneous harmonic polynomial $P(z)$ to be HN 
can be obtained or described 
in terms of the graph structure 
of $\cG(P)$. 

We first give in Subsection \ref{S8.1} the definition 
of the graph $\cG(P)$ for any homogeneous  
harmonic polynomial $P(z)$ and discuss 
the connectedness 
reduction (see Corollary \ref{conn-reduction}), 
i.e. a reduction of the VC 
to the homogeneous HN polynomials $P$ 
such that $\cG(P)$ is connected.  
We then consider in Subsection \ref{S8.2} 
a connection of $\cG(P)$ 
with the tree expansion formula 
derived in \cite{Me} and \cite{TreeExp} for the inversion 
pair $Q(z)$ of $P(z)$ (see Proposition \ref{grf-MainThm-2}). 
As an application of the connection, we give another 
proof for the connectedness 
reduction given in Corollary \ref{conn-reduction}.

\subsection{Definition and the Connectedness Reduction}
\label{S8.1}

For any $\beta \in \bC^n$, 
set $h_\beta(z)\!:=\la \beta, z\ra$ and 
$\beta_D\!:=\la \beta, D\ra$, where 
$\la\cdot, \cdot\ra$ is the standard 
$\bC$-bilinear form of $\bC^n$. Let $X(\bC)$ 
denote the set of all isotropic elements 
of $\bC^n$, i.e. the set of all elements 
$\alpha \in \bC^n$ such that 
$\la \alpha, \alpha \ra=0$. 

Recall that we have the following 
fundamental theorem on homogeneous  
harmonic polynomials. 

\begin{theo}\label{T8.1.1}
For any homogeneous harmonic polynomial 
$P(z)$ of degree $d\geq 2$, we have
\begin{align} \label{d-Form-0}
P(z)=\sum_{i=1}^k c_i h_{\alpha_i}^d (z)
\end{align}
for some $c_i\in \bC^\times$ and 
$\alpha_i\in X(\bC^n)$ $(1\leq i\leq k)$.
\end{theo}

Note that, replacing $\alpha_i$ in Eq.\,(\ref{d-Form-0})
by $c_i^{-\frac 1d}\alpha_i$, we may also write $P(z)$ as
\begin{align} \label{d-Form}
P(z)=\sum_{i=1}^k h_{\alpha_i}^d (z)
\end{align}
with $\alpha_i\in X(\bC^n)$ $(1\leq i\leq k)$.

For the proof of Theorem \ref{T8.1.1}, 
see, for example, \cite{I} 
and \cite{Roel}.

We fix a homogeneous harmonic polynomial 
$P(z)\in \bCz$ of degree $d\geq 2$, and 
assume that $P(z)$ is given by Eq.\,(\ref{d-Form}) 
for some $\alpha_i \in X(\bC^n)$ $(1\leq i\leq k)$. 
We may and will always assume 
$\{h_{\alpha_i}^d(z) | 1\leq i\leq k\}$
are linearly independent in $\bCz$.

Recall the following matrices had been introduced 
in \cite{HNP}:
\begin{align}
A_P&=(\la \alpha_i, \alpha_j\ra)_{k\times k}, \\
\Psi_P&=(\la \alpha_i, \alpha_j \ra h_{\alpha_j}^{d-2}(z) )_{k\times k}.
\end{align}

Then we have the following criterion 
of Hessian nilpotency for homogeneous 
harmonic polynomials. For its proof, 
see Theorem $4.3$ in \cite{HNP}.

\begin{theo}\label{Crit-2}
Let $P(z)$ be as above. 
Then, for any $m\geq 1$, 
we have 
\begin{align}\label{E6.1.7}
\text{Tr\,} \Hes^m (P)=(d(d-1))^m \text{Tr\,} \Psi_P^m.
\end{align}
In particular, $P(z)$ is HN if and only if 
the matrix $\Psi_P$ is nilpotent.
\end{theo}

One simple remark on the criterion above is as follows.

Let $B$ be the $k\times k$ diagonal matrix with 
the $i^{th}$ $(1\leq i\leq k)$ diagonal entry being 
$h_{\alpha_i}(z)$. For any $1\leq j\leq k$, set
\begin{align}
\Psi_{P; j}\!:= B^j A_P B^{d-2-j}=
(h_{\alpha_i}^j \la \alpha_i, \alpha_j \ra h_{\alpha_j}^{d-2-j}).
\end{align}

Then, by repeatedly applying the fact that, 
for any two $k\times k$ matrices $C$ and $D$, 
$CD$ is nilpotent iff so is $DC$, it is easy to
see that Theorem  \ref{Crit-2} can also 
be re-stated as follows.

\begin{corol}
Let $P(z)$ be given by Eq.\,$(\ref{d-Form})$ with $d\geq 2$. 
Then, for any $1\leq j\leq d-2$ and $m\geq 1$, 
we have 
\begin{align}
\text{Tr\,} \Hes^m (P)=(d(d-1))^m \text{Tr\,} \Psi_{P; j}^m.
\end{align}
In particular, $P(z)$ is HN if and only if 
the matrix $\Psi_{P; j}$ is nilpotent.
\end{corol}

Note that, when $d$ is even, we may choose $j=(d-2)/2$. 
So $P$ is HN iff the symmetric matrix 
\begin{align}\label{(d-2)}
\Psi_{P; (d-2)/2}(z)=(\, h_{\alpha_i}^{(d-2)/2}(z)\, 
\la \alpha_i, \alpha_j \ra \,
h_{\alpha_j}^{(d-2)/2}(z)\,)
\end{align}
is nilpotent.

Motivated by the criterion above, 
we assign a graph $\cG(P)$ 
to any homogeneous harmonic 
polynomial $P(z)$ as follows.

We fix an expression as in 
Eq.\,$(\ref{d-Form})$ for $P(z)$. 
The set of vertices of $\cG(P)$ 
will be the set of positive integers 
$[{\bf k}]\!:=\{1, 2, \dots, k\}$.
The vertices i and j of $\cG(P)$ are 
connected by an edge iff 
$\la \alpha_i, \alpha_j \ra \neq 0$. 
In this case, we get a finite graph. 

Furthermore, we may also label edges of 
$\cG(P)$ by assigning 
$\la \alpha_i, \alpha_j \ra$ or 
$(h_{\alpha_i}^{(d-2)/2} \la \alpha_i, 
\alpha_j \ra h_{\alpha_i}^{(d-2)/2})$,
when $d$ is even, for the edge 
connecting vertices $i, j\in [\bf k]$. 
We then get a labeled graph whose 
adjacency matrix is exactly $A_P$ 
or $\Psi_{P, (d-2)/2}$ (depending on the labels we choose 
for the edges of $\cG(P)$). 

Naturally, one may also ask the following (open) questions.

\begin{prob}\label{grf-prob}
$(a)$ Find some necessary or sufficient conditions on 
the $($labeled$)$ graph $\cG(P)$ such that 
the homogeneous harmonic 
polynomial $P(z)$ is HN.

$(b)$ Find some necessary or sufficient conditions 
on the $(labeled)$ graph $\cG(P)$ such that 
the VC holds for the homogeneous 
HN polynomial $P(z)$.
\end{prob}

First, let us point out that, to approach 
the open problems above, it will be enough 
to focus on homogeneous harmonic polynomials $P$ 
such that the graph $\cG(P)$ is connected. 
  
Suppose that the graph $\cG(P)$ 
is a disconnected graph with $r\geq 2$ 
connected components. 
Let $[{\bf k}]=\sqcup_{i=1}^r I_i$ be the corresponding 
partition of the set 
$[\bf k]$ of vertices of $\cG(P)$. 
For each $1\leq i \leq r$, we set 
$P_i(z)\!:= \sum_{\alpha \in I_i} 
h_\alpha^d(z)$. 

Note that, by Lemma \ref{p-inv-L4}, 
$P_i$ $(1\leq i\leq r)$ are disjoint 
to each other, so Corollary \ref{p-inv-MainCorol} 
applies to the sum $P=\sum_{i=1}^r P_i$. 
In particular, we have, 

{\it
$(a)$ $P$ is HN iff each $P_i$ is HN. 

$(b)$ if the VC holds for each $P_i$, 
then it also holds for $P$.}
\vskip2mm
Therefore, we have the following {\it connectedness reduction}.

\begin{corol}\label{conn-reduction}
To study homogeneous HN polynomials $P$ 
or the VC for homogeneous HN polynomials $P$, it will 
be enough to consider the case 
when $\cG(P)$ is connected.
\end{corol}

Note that, the property $(a)$ above 
was first proved by R. Willems  
(\cite{Roel}) by using the criterion in 
Theorem  \ref{Crit-2}. 
$(b)$ was first proved 
by the author by a different argument, and 
with the author's permission, 
it had also been included in \cite{Roel}. 

Finally, let us point out that R. Willems  (\cite{Roel}) 
has proved the following very interesting results 
on Open Problem \ref{grf-prob}.  

\begin{theo}$($\cite{Roel}$)$
Let $P$ be a homogeneous HN polynomial as 
in Eq.$(\ref{d-Form})$ with $d\geq 4$. 
Let $l(P)$ be the dimension of 
the vector subspace of $\bC^n$ 
spanned by $\{\alpha_i\,|\, 1\leq i\leq k\}$.
Then 
\begin{enumerate}
\item If $l(P)=1, 2, k-1$ or $k$, the graph 
$\bG(P)$ is totally disconnected 
$($i.e. $\cG(P)$ is the graph with no edges$)$.

\item If $l(P)=k-2$ and $\cG(P)$ is connected, 
then $\cG(P)$ is the complete bi-graph $K(4, k-4)$.

\item In the case of $(a)$ and $(b)$ above, 
the VC holds.
\end{enumerate}
\end{theo}

Furthermore, it has also been shown in \cite{Roel} that, 
for any homogeneous HN polynomials $P$, the graph
$\cG(P)$ can not be any path nor cycles 
of any positive length. 
For more details, see \cite{Roel}.

\subsection{Connection with the Tree Expansion Formula 
of Inversion Pairs}\label{S8.2}

First let us recall the tree expansion formula 
derived in \cite{Me}, \cite{TreeExp} for 
the inversion pair $Q(z)$.

Let $\bT$ denote the set of all trees, i.e. 
the set of all connected and simply connected 
finite simple graphs. For each tree 
$T\in \bT$, denote by $V(T)$ and $E(T)$ 
the sets of all vertices and edges of $T$, 
respectively. Then we have the following 
tree expansion formula 
for inversion pairs.

\begin{theo} $($\cite{Me}, \cite{TreeExp}$)$ 
Let $P\in \bCzz$ with $o(P)\geq 2$ and $Q$ its inversion pair. 
For any $T\in \bT$, set 
\begin{align}
Q_{T, P}=\sum_{\ell: E(T)\to [\bf n]} 
\prod_{v\in V(T)} D_{adj(v), \ell} P,
\end{align}
where $adj(v)$ is the set 
$\{e_1, e_2,\dots, e_s\}$ of edges of $T$ 
adjacent to $v$, and 
$D_{adj(v), \ell}=D_{\ell(e_1)}D_{\ell(e_2)} \cdots D_{\ell(e_s)}$.

Then  the inversion pair $Q$ of $P$ is given by
\begin{align}\label{Me-Wr-eq1}
Q=\sum_{T\in \bT} \frac 1{|Aut(T)|} Q_{T, P}.
\end{align}
\end{theo}

Now we assume $P(z)$ is a homogeneous harmonic polynomial 
$d\geq 2$ and has expression in Eq.\,$(\ref{d-Form})$. 
Under this assumption, it is easy to 
see that $Q_{T, P}$ $(T\in \bT)$ becomes
\begin{align}\label{QTP-1}
Q_{T, P}=\sum_{f: V(T)\to [\bf k]}
\sum_{\ell: E(T)\to [\bf n]} \prod_{v\in V(T)} 
D_{adj(v), \ell} h_{\alpha_{f(v)}}^d(z).
\end{align}

The role played by the graph $\cG(P)$ 
of $P$ is to restrict the maps 
$f: V(T)\to V(\cG(P))(=[\bf k])$ in 
Eq.\,$(\ref{QTP-1})$ to a special 
family of maps. To be more precise, 
let $\Omega(T, \cG(P))$ be the set 
of maps $f: V(T) \to [\bf k]$ 
such that, {\it for any distinct adjoint 
vertices $u, v\in V(T)$, 
$f(u)$ and $f(v)$ are distinct and 
adjoint in $\cG(P)$}. 
Then we have the following lemma.

\begin{lemma}\label{grp-8.2-L1}
For any $f: V(T) \to [\bf k]$ with 
$f \not \in \Omega(T, \cG(P))$, we have
\begin{align}
\sum_{\ell: E(T)\to [\bf n]} \prod_{v\in V(T)} 
D_{adj(v), \ell} h_{\alpha_{f(v)}}^d(z)=0.
\end{align}
\end{lemma}

\pf Let $f:V(T) \to [\bf k]$ as in the lemma. 
Since $f \not \in \Omega(T, \cG(P))$, there exist distinct
adjoint $v_1, v_2\in V(T)$ such that, either 
$f(v_1)=f(v_2)$ or $f(v_1)$ and $f(v_2)$ 
are not adjoint in the graph $\cG(P)$. 
In any case, we have 
$\la \alpha_{f(v_1)}, \alpha_{f(v_2)} \ra =0$.

Next we consider contributions to the RHS of 
Eq.\,$(\ref{QTP-1})$ 
from the vertices $v_1$ and $v_2$. 
Denote by $e$ the edge of $T$ connecting 
$v_1$ and $v_2$, 
and $\{e_1, \dots e_r\}$ 
(resp.\,$\{\tilde e_1, \dots \tilde e_s\}$) 
the set of edges connected with $v_1$ 
(resp.\,$v_2$) besides the edge $e$.
Then, for any $\ell: E(T)\to [\bf n]$, 
the factor in the RHS of Eq.\,$(\ref{QTP-1})$ 
from the vertices $v_1$ and $v_2$ 
is the product
\begin{align}\label{PFEQ-1}
\left (D_{\ell(e)} D_{\ell(e_1)} \cdots D_{\ell(e_r)} 
h_{\alpha_{f(v_1)}}^d(z)\right) \left( D_{\ell(e)} 
D_{\ell(\tilde e_1)} \cdots 
D_{\ell(\tilde e_s)} h_{\alpha_{f(v_2)}}^d(z)\right).
\end{align}

Define an equivalent relation 
for maps $\ell: E(T) \to [\bf n]$ 
by setting $\ell_1 \sim \ell_2$ iff 
$\ell_1$, $\ell_2$ have same image
at each edge of $T$ except $e$.
Then, by taking sum of the terms in 
Eq.\,$(\ref{PFEQ-1})$ over each equivalent class, 
we get the factor 
\begin{align}\label{PFEQ-2}
\left \la \nabla  D_{\ell(e_1)} \cdots D_{\ell(e_r)} 
h_{\alpha_{f(v_1)}}^d(z), \,
\nabla D_{\ell(\tilde e_1)} \cdots 
D_{\ell(\tilde e_s)} h_{\alpha_{f(v_2)}}^d(z)\right \ra.
\end{align}

Note that $D_{\ell(e_1)} \cdots D_{\ell(e_r)} 
h_{\alpha_{f(v_1)}}^d(z)$ and $D_{\ell(\tilde e_1)} \cdots 
D_{\ell(\tilde e_s)} h_{\alpha_{f(v_2)}}^d(z)$ 
are constant  multiples of some integral powers of 
$h_{\alpha_{f(v_1)}}(z)$ and $h_{\alpha_{f(v_2)}}(z)$, 
respectively. Therefore, 
$\la \alpha_{f(v_1)}, \alpha_{f(v_2)}\ra(=0)$
appears as a multiplicative constant factor 
in the term in Eq.\,$(\ref{PFEQ-2})$,
which makes the term zero. 
Hence the lemma follows.
\epfv

One immediate consequence of the lemma above 
is the following proposition.

\begin{propo}\label{grf-MainThm-2}
With the setting and notation as above, we have
\begin{align}\label{QTP-2}
Q_{T, P}=\sum_{f\in \Omega(T, \cG(P))}\,
\sum_{\ell: E(T)\to [\bf n]} \prod_{v\in V(T)} 
D_{adj(v), \ell} h_{\alpha_{f(v)}}^d(z).
\end{align}
\end{propo}

\begin{rmk}
$(a)$ For any $f\in \Omega(T, \cG(P))$, 
$\{ f^{-1}(j)\,|\, j \in \mbox{Im}(f) \}$ 
gives a partition of $V(T)$ since 
no two distinct vertices in $f^{-1}(j)$ 
$(j \in \mbox{Im}(f))$ can be adjoint. 
In other words, $f$ is nothing 
but a proper coloring for the tree $T$,
which is also subject to certain 
more conditions from the graph 
structure of $\cG(P)$. 
It is interesting to see that
the coloring problem of graphs 
also plays a role in the inversion problem 
of symmetric formal maps. 

$(b)$ It will be interesting to see if more results can 
be derived from the graph $\cG(P)$ via the formulas in
Eqs.\,$(\ref{Me-Wr-eq1})$ and $(\ref{QTP-2})$. 
\end{rmk}

\begin{rmk}\label{PST-2}
By similar arguments 
as those in proofs of Lemma \ref{grp-8.2-L1}, 
one may get another proof for 
Theorem \ref{p-inv-MainThm-2}  
in the setting as in Lemma \ref{p-inv-L4}.
\end{rmk}

Finally, as an application of Proposition \ref{grf-MainThm-2} above, 
we give another proof for the connectedness reduction given 
in Corollary \ref{conn-reduction}.

Let $P$ as given in Eq.\,$(\ref{d-Form})$ 
with the inversion pair $Q$. 
Suppose that there exists 
a partition $[{\bf k}]=I_1 \sqcup I_2$ with 
$I_i \neq \emptyset$. 
Let $P_i=\sum_{\alpha\in I_i} h_{\alpha}^d(z)$
$(i=1,2)$ and $Q_i$ the inversion pair of $P_i$. 
Then we have $P=P_1+P_2$ and 
$\cG(P_1)\sqcup \cG(P_2)=\cG(P)$.
Therefore, to show the connectedness reduction 
discussed in the previous subsection, 
it will be enough to show $Q=Q_1+Q_2$.
But this will follow immediately from 
Eqs.\,$(\ref{Me-Wr-eq1})$, $(\ref{QTP-2})$ 
and the following lemma.

\begin{lemma}
Let $P$, $P_1$ and $P_2$ as above, 
then, for any tree $T\in \bT$, we have
\begin{align*}
\Omega(T, \cG(P))=\Omega(T, \cG(P_1)) \sqcup \Omega(T, \cG(P_2)).
\end{align*}
\end{lemma}
\pf For any $f\in \Omega(T, \cG(P))$, 
$f$ preserves the adjacency of vertices 
of $\cG(P)$. Since $T$ as a graph 
is connected, 
$\mbox{Im}(f)\subset V(\cG(P))$ as 
a (full) subgraph of $\cG(P)$ must 
also be connected. 
Therefore,  $\mbox{Im}(f)\subset V(\cG(P_1))$ or
$\mbox{Im}(f)\subset V(\cG(P_2))$. 
Hence $\Omega(T, \cG(P))\subset \Omega(T, \cG(P_1)) \sqcup \Omega(T, \cG(P_2))$. 
The other way of containess is obvious.
\epfv

{\small \sc Department of Mathematics, Illinois State University,
Normal, IL 61790-4520.}

{\em E-mail}: wzhao@ilstu.edu.

\end{document}